\documentclass[12pt,draftcls,onecolumn]{IEEEtran}
\usepackage{amsfonts,amssymb,epsfig, float}
\usepackage{times}                                     % if you need a4paper
\usepackage[centertags]{amsmath}
                      \usepackage{amsfonts,amssymb}                                    % needed if you want to
at 17pt

\catcode`\@=11
\def\downparenfill{$\m@th\braceld\leaders\vrule\hfill\bracerd$}
\def\overparen#1{\mathop{\vbox{\ialign{##\crcr\crcr \noalign{\kern0.4ex}
\downparenfill\crcr\noalign{\kern0.4ex\nointerlineskip}
$\hfil\displaystyle{#1}\hfil$\crcr}}}\limits} \catcode`\@=12
\catcode`\@=11
\let\ps@plainbuf\ps@plain
\def\ps@plain{\ps@plainbuf
              \ifnum \value{page}=1
                  \def\@oddfoot{\@foottitle}
              \fi}
\def\foottitle#1{\gdef\@foottitle{\hbox\sl\hfill#1\hfill}}
\gdef\@foottitle{\rm\hfil\thepage\hfil} \catcode`\@=12
\newtheorem{theorem}{Theorem}

\newtheorem{lemma}[theorem]{Lemma}
\newtheorem{corollary}[theorem]{Corollary}

\newtheorem{assumption}[enumi]{Assumption}
\newtheorem{remark}[theorem]{Remark}

\begin{document}
\centerline{\Large\bf Further Constructions of Control-Lyapunov
Functions}\smallskip\centerline{\Large\bf and Stabilizing
Feedbacks for Systems Satisfying}\smallskip\centerline{\Large\bf
the Jurdjevic-Quinn Conditions
%Control-Lyapunov Functions for Systems Satisfying
%the}\smallskip\centerline{\Large\bf Conditions of the
%Jurdjevic-Quinn Theorem
}
\medskip\medskip\medskip
\centerline{\large Fr\'{e}d\'{e}ric Mazenc and Michael Malisoff
 \footnote{\  Corresponding Author: F. Mazenc.  The first author  was supported by the MERE Project.
The second author was supported by NSF Grant 0424011. F. Mazenc is
with the Projet MERE INRIA-INRA, UMR Analyse des Syst\`{e}mes et
Biom\'{e}trie INRA, 2, pl. Viala,
        34060 Montpellier, France,
        $\mathtt{mazenc@helios.ensam.inra.fr}$.
M. Malisoff is with the Department of Mathematics, Louisiana State
University, Baton Rouge, LA 70803-4918,
$\mathtt{malisoff@lsu.edu}$. }}
\medskip\medskip\medskip\medskip\medskip

%%%%%%%%%%%%%%%%%%%%%%%%%%%%%%%%%%%%%%%%%%%%%%%%%%%%%%%%%%%%%%%%%%%%%%%%%%%%%%%%

{\abstract For a broad class of nonlinear systems, we construct
smooth control-Lyapunov functions whose derivatives along the
trajectories of the systems can be made negative definite by
smooth control laws that are arbitrarily small in norm.  We assume
our systems satisfy appropriate generalizations of the
Jurdjevic-Quinn conditions. We also design state feedbacks of
arbitrarily small norm that render our systems
integral-input-to-state stable to actuator errors.
\medskip

}

\medskip\medskip\medskip
%%%%%%%%%%%%%%%%%%%%%%%%%%%%%%%%%%%%%%%%%%%%%%%%%%%%%%%%%%%%%%%%%%%%%%%%%%%%%%%%
\section{Introduction}

Lyapunov stability is of paramount importance in
nonlinear control theory. In many important applications, it is
very beneficial to have  a continuously differentiable Lyapunov
function whose derivative along the trajectories of the system can
be made negative definite by an appropriate choice of feedback.
Observe in particular that:
\newline
\noindent $\bullet$ Recent advances in the stabilization of
nonlinear delay systems (e.g., \cite{Jank2,MMN,Teel:98})
are based on knowledge of continuously differentiable Lyapunov
functions.
\newline
\noindent $\bullet$  Lyapunov functions are  very efficient tools
for robustness analysis. For example, many proofs of nonlinear
disturbance-to-state $L^p$ stability properties rely on Lyapunov
functions; see \cite[Chapter 13]{Isi} and \cite{ASW00, Lichi,MazP}.
Moreover,  control-Lyapunov function (CLF) based control
designs guarantee robustness to different types of deterministic
\cite{FK} and stochastic disturbances, and to unmodeled dynamics
\cite{PraW,SepJK}.
\newline \noindent $\bullet$ When a CLF satisfying the {\em small control property} (as defined below)
is available,  the universal formula  in \cite{Sontag} provides an
explicit expression for a stabilizing feedback that is also an
optimal control for a  suitable optimization problem whose value
function is the CLF;  see \cite{Sontag}.
\newline \noindent $\bullet$
Backstepping and  forwarding require Lyapunov functions of class
$C^1$ for the subsystems \cite{SepJK}.

The converse Lyapunov theorem (see \cite{Kur}) ensures that, for any
system that is globally asymptotically stabilizable by $C^1$
feedback, a CLF exists. Unfortunately, for nonlinear control
systems, determining {\em explicit expressions} for CLFs is in
general difficult. Fortunately, for large classes of systems, one
can determine functions whose derivatives along the trajectories can
be rendered negative {\em semi}-definite. If the systems satisfy the
so-called weak Jurdjevic-Quinn conditions (defined below),
which generalize those given in \cite{Jurd}, then globally
asymptotically stabilizing feedbacks can be constructed.
However, in this case, explicit formulas for CLFs are not
generally available. This motivates the following fundamental
question: {\em When the Jurdjevic-Quinn method applies, is it
possible to design explicit CLFs?}

In \cite{FauPo}, where this issue was addressed for the first time,
a method was presented for designing explicit CLFs for affine
homogeneous systems that satisfy the Jurdjevic-Quinn conditions. Our
objective in the present note is to extend the main result of
\cite{FauPo} by constructing CLFs for systems satisfying appropriate
generalizations of the Jurdjevic-Quinn conditions, but not
necessarily having the homogeneity property, {}including cases where
the system may not be control-affine. Our work also complements
\cite{MazD} where strong Lyapunov functions are constructed for a
large family of systems satisfying either an appropriate Lie
algebraic condition or which can be shown to be stable using the
LaSalle invariance principle. The main difference between the
present work and \cite{MazD} is that in \cite{MazD}, only systems
without input are considered whereas here we consider systems with
input.

We end this introduction by recalling some basic facts on the
Jurdjevic-Quinn method. We say {}(see for example \cite{FauPo} for
the relevant definitions) that a nonlinear control-affine system
\begin{equation}
\label{1gh}
\dot{x} \; = \; f(x) + g(x) u \; , \; g(x) = (g_1(x), \dots ,g_m(x))
\end{equation}
satisfies the {}{\em (weak) Jurdjevic-Quinn conditions} provided
there exists a function  $V: {\mathbb R}^n \rightarrow {\mathbb R}$
satisfying the following three properties: $V$ is positive definite
and radially unbounded; for all $x \in {\mathbb R}^n$, $L_f V(x)
\leq 0$; and there exists an integer $l$ such that the set
\[
W(V) = \left\{x \in {\mathbb R}^n : \forall k \in\{ 1, \dots, m\}
{\rm \ and\ } i \in  \{0, \dots ,l\},  L_f V(x) =
L_{ad_{f}^i(g_k)}V(x) = 0\right\}
\]
equals $\{0\}$. {}Here and in the sequel, we assume all functions
are sufficiently smooth.  If (\ref{1gh}) satisfies the weak
Jurdjevic-Quinn conditions, then it is globally asymptotically
stabilized by any feedback $u = - \xi(x) L_g V(x)^\top$ where $\xi$
is any positive function of class $C^1$. The proof of this result
relies on the LaSalle Invariance Principle.

The remainder of this paper is organized as follows.  In Section
\ref{sec1}, we present our main result. Section \ref{discussion} is
devoted to a discussion of our main result, Section \ref{proof} to
its proof, and Section \ref{secex} to an illustrating example.
Section \ref{seciss} constructs feedbacks for our systems that have
arbitrarily small norm  and that in addition achieve
integral-input-to-state stability relative to actuator errors.
Concluding remarks in Section {}\ref{seconc} end our work.

\section{Main result}
\label{sec1}

Recall (cf. \cite{AAS02}) that a $C^1$ positive definite function
$V(\cdot)$ on ${\mathbb R}^n$ is called a {\em control-Lyapunov
function (CLF)} for a system $\dot{\chi}  =  \varphi_1(\chi) +
\varphi_2(\chi)u$ with input $u$ provided it is radially unbounded
and satisfies $ L_{\varphi_1} V(\chi) \geq 0 \; \Rightarrow \;
\left[\chi = 0 \quad \mbox{or} \quad L_{\varphi_2} V(\chi) \neq
0\right] $. We use $\dot V(x,u)$ to denote the derivative $\dot
V(x,u)=L_{\varphi_1}V(x)+L_{\varphi_2}V(x)u$  of $V$ along
trajectories of the system. We often suppress the arguments of $\dot
V$ to simplify the notation. We say that a CLF $V(\cdot)$ for the
system $\dot{\chi}  =  \varphi_1(\chi) + \varphi_2(\chi)u$ satisfies
the {\em small control property} \cite{Sontag} provided for each
$\varepsilon > 0$, there exists $\delta(\varepsilon) > 0$ such that
if $0 < |\chi| < \delta(\varepsilon)$, then there exists $u$
(possibly depending on $\chi$) such that $|u| < \varepsilon$ and
$L_{\varphi_1} V(\chi) + L_{\varphi_2} V(\chi)u < 0$.

We next provide our main CLF and stabilizing feedback constructions
for the fully nonlinear system
\begin{equation}
\label{fn}
\dot x = F(x,u)
\end{equation}
where $x\in {\mathbb R}^n$, $u\in {\mathbb R}^m$ is the control,
$F(0,0)=0$,  and the function $F$ is assumed to be $C^1$.  We
further assume that $u\mapsto F(x,u)$ is $C^2$ (i.e., the second
order partial derivatives, with respect to the components of $u$,
of each component of $F$ are continuous), so the functions
\begin{equation} \label{data} f(x):=F(x,0),\; \; \; \;
g(x):=\frac{\partial F}{\partial u}(x,0)
\end{equation}
are at least $C^1$.  Finally, we assume:

\begin{assumption}
\label{H1} A smooth function $V(x)$ that is radially unbounded and
positive definite and such that
\begin{equation}
\label{2gh} L_{f} V(x)  \; \leq \; 0 \; \; \; \forall x\in
{\mathbb R}^n
\end{equation}
is known.\end{assumption}

\begin{assumption}
\label{H2}
A vector field $G(x)$ such that if
$L_{g}V(x) = 0$ and $x \neq 0$, then we either have
$L_{f}L_{G}V(x) < 0$ or $L_{f}V(x) < 0$, is known.
\end{assumption}

We are ready to state our main result.
\begin{theorem}
\label{thgh} Assume  the data (\ref{data}) satisfy Assumptions
\ref{H1}-\ref{H2}. Then one can determine a positive definite
smooth function $\delta:[0,\infty)\to[0,\infty)$ and a function
$\Omega:[0,\infty)\to [0,\infty)$ such that
\begin{equation}
\label{3gh} V^\sharp(x) \; = \; V(x) + \int_0^{V(x)}\Omega(s){\rm
d}s+\delta(V(x)) L_{G}V(x)
\end{equation}
is a CLF for  (\ref{fn}) that satisfies the small control property.
In fact, for each real-valued $C^1$ positive function
$\bar\xi(\cdot)$, one can determine a function $\delta(\cdot)$, and
a $C^1$ function $\xi:[0,\infty)\to (0,\infty)$ satisfying
$\xi(s)\le \bar \xi(s)$ for all $s \ge 0$, such that (\ref{3gh}) is
a CLF for (\ref{fn}) satisfying the small control property
 whose derivative along
the trajectories of (\ref{fn}) in closed-loop with the feedback
\begin{equation}
\label{bj1} u \; = \; - \xi(V(x)) L_{g}V(x)^\top
\end{equation}
is negative definite.
%In particular, (\ref{fn}) admits smooth
%globally  stabilizing state feedbacks of arbitrarily small
%amplitude.
\end{theorem}

\section{Discussion of Theorem \ref{thgh}}
\label{discussion}

\noindent 1. Assumptions \ref{H1} and \ref{H2} are similar to the
assumptions of the main result of \cite{FauPo}. In particular, for
the special case where $F$ is control-affine, \cite{FauPo}
provides an explicit expression for a vector field $G(x)$ such
that Assumption \ref{H2} holds whenever the so-called
``weak Jurdjevic-Quinn conditions" (see the introduction) are satisfied.
This vector field is
not continuous at the origin but it turns out that there exists an
integer $N \geq 1$ such that the vector field $G_N(x) = V(x)^N
G(x)$ is of class $C^{\infty}$ for $V$ satisfying our assumptions.
The equality $L_f L_{G_N} V(x)  =  N V(x)^{N - 1} L_f V(x) L_G
V(x)
  +  V(x)^N L_f L_{G} V(x)$ then implies that if $G(x)$
satisfies Assumption \ref{H2}, and if Assumption \ref{H1} also
holds, then $G_N(x)$ satisfies Assumption \ref{H2} as well.
Consequently, one can take advantage of the formula in
\cite{FauPo} to determine a $C^\infty$ vector field  for which
Assumption \ref{H2} is satisfied.

\noindent 2. No restriction on the size of the function
$\xi(\cdot)$ in (\ref{bj1}) is imposed. Therefore, the family of
feedbacks (\ref{bj1}) contains elements that are arbitrarily small
in (sup) norm.  In fact, for any continuous positive function
$\epsilon:[0,\infty)\to (0,\infty)$, we can design our stabilizing
feedback $u$ so that it satisfies $|u(x)|\le \epsilon(|x|)$ for
all $x\in {\mathbb R}^n$.

\noindent 3. An important class of dynamics covered by Theorem
\ref{thgh} is described by the so-called {\em Euler-Lagrange
equations}
\begin{equation}
\label{el} \frac{d}{dt}\left(\frac{\partial L}{\partial \dot
q}(q,\dot q)\right)-\frac{\partial L}{\partial q}(q,\dot q)=\tau
\end{equation}
for the motion of mechanical systems,  in which $q$ represents the
generalized configuration coordinates, $L=K-P$ is the difference
between the kinetic energy $K$ and potential energy $P$, and
$\tau$ is the control \cite{V00}.  In  standard cases,  $K(q,\dot
q)=\frac{1}{2}\dot q^\top M(q)\dot q$ where the inertia matrix
$M(q)$ is $C^1$ and everywhere symmetric and positive definite.
Then the generalized momenta $\frac{\partial L}{\partial \dot q}$
are given by $p=M(q)\dot q$, so in terms of the state  $x=(q,p)$,
the equations (\ref{el}) become \cite{V00}
\begin{equation}
\label{hamil} \dot q=\frac{\partial H}{\partial
p}(q,p)^\top=M^{-1}(q)p, \; \; \; \; \dot p=-\frac{\partial
H}{\partial q}(q,p)^\top+\tau,
\end{equation}
where $H(q,p)=\frac{1}{2}p^\top M^{-1}(q)p+P(q)$ is the total energy
of the system.  We make the following additional assumptions: (a)
$P(q)$ is positive definite and radially unbounded and (b) $\nabla
P(q)\ne 0$ whenever $q\ne 0$.   (These two assumptions are not too
restrictive since one can often modify $H$ and $\tau$ to get a new
system that satisfies these assumptions.  Condition (a) can be
weakened by assuming there is a constant $c$ such that
$q \mapsto P(q)+c$ is radially unbounded and positive definite
in which case we simply add $c$ to the function $V$ in what follows.)
Then $H$ is positive definite and radially unbounded, so  $V = H$
satisfies Assumption \ref{H1}. The radial unboundedness follows from
the continuity of the (positive) eigenvalues of the positive
definite matrix $M^{-1}(q)$ as functions of
$q$ \cite[Appendix A4]{S98a}, which implies that each compact set
$\cal S$ of $q$ values admits a
constant $c_{\cal S}>0$ such that $p^\top M^{-1}(q)p\ge c_{\cal
S}|p|^2$ for all $q\in {\cal S}$ and all $p$. In our general
notation with $x=(q,p)$, we get $L_fV(x)\equiv 0$ and
$L_gV(x)=H_p(x)=p^\top M^{-1}(q)$. Choosing $ G(x) = [\; 0\; \;
\nabla P(q)^\top\; ]^\top$ gives $L_GV(x)=H_p(x)\nabla P(q)^\top$.
Therefore,  if $L_gV(x)=p^\top M^{-1}(q)= 0$ and $x\ne 0$, then
$p=0$
 and therefore also $L_fL_GV(x)=-\nabla P(q) M^{-1}(q)\nabla
P(q)^\top$ and $q\ne 0$. Since $M^{-1}$ is everywhere positive
definite, Assumption \ref{H2} therefore reduces to our assumption
(b) and therefore is satisfied as well. We study a special case of
(\ref{hamil}) in Section \ref{secex} below, where we explicitly
compute the corresponding CLF  and stabilizing feedback.

\section{Proof of Theorem \ref{thgh}}
\label{proof} \subsection{Control Affine Case}\label{ca} We fix a
positive function $\bar \xi:[0,\infty)\to(0,\infty)$, and
functions $V$ and $G$ satisfying Assumptions \ref{H1}-\ref{H2}. We
begin by proving Theorem \ref{thgh} for the case where (\ref{fn})
is control affine, i.e., of the form (\ref{1gh}).
In this control affine case, the conclusions of our theorem will
hold with $\Omega\equiv 0$ and $\xi\equiv \bar \xi$.  In Section
\ref{fnc}, we will modify our constructions to handle the fully
nonlinear system (\ref{fn}).

\noindent {\em First step.} We exhibit a family of functions
$\delta(\cdot)$ for which the function
\begin{equation}
\label{3gha} U(x) \; := \; V(x) +\delta(V(x)) L_{G}V(x)
\end{equation} is positive definite and radially unbounded.
One can determine $\alpha_i(\cdot)$ of class ${\cal K}_{\infty}$
such that $\alpha_{1}(|x|) \leq V(x) \leq \alpha_{2}(|x|)$ and
$|L_{G}V(x)| \leq \alpha_{3}(|x|)$ for all $x\in {\mathbb R}^n$.
It follows that
\begin{equation}
\label{6gh} U(x) \;  \geq \;  \alpha_{1}(|x|) -
\delta(V(x))\alpha_3(|x|)\; \ge \; \alpha_{1}(\alpha_{2}^{-
1}(V(x)))\; - \; \delta(V(x))\alpha_{3}(\alpha_{1}^{- 1}(V(x)))
\end{equation}
for all $x\in {\mathbb R}^n$. We can use standard results to find
a $C^1$ function $\delta:[0,\infty)\to[0,\infty)$  such that
\begin{equation}
\label{7gh} \delta(v) \leq \frac{\alpha_{1}(\alpha_{2}^{-
1}(v))}{1+2 \alpha_{3}(\alpha_{1}^{- 1}(v))} \quad   \forall v \ge
0 \ .
\end{equation}
With such a  function $\delta(\cdot)$, the inequality $U(x) \geq
 \frac{1}{2}\alpha_{1}(\alpha_{2}^{- 1}(V(x)))$ for all $x\in
{\mathbb R}^n$ follows from (\ref{6gh}). Since $V(x)$ is positive
definite and radially unbounded and $
\frac{1}{2}\alpha_{1}(\alpha_{2}^{- 1}(\cdot))$ is of class ${\cal
K}_{\infty}$, this implies that $U(x)$ is positive definite and
radially unbounded as well. In the next steps, we impose further
restrictions on $\delta$.

\noindent {\em Second step.} Along the trajectories $x(t)$ of our
system (\ref{1gh}) in closed-loop with the feedback
$u=-\bar\xi(V(x))L_gV(x)^\top$, the derivative $\dot U$ of $U(x)$
from (\ref{3gha}) reads \begin{eqnarray} \label{3h0} \dot U  &=&
\left[L_{f}V(x) - \bar\xi(V(x))
|L_{g}V(x)|^2\right]\left[1 + \delta'(V(x)) L_{G}V(x)\right]\nonumber\\
&&+ \; \delta(V(x)) L_{f} L_{G}V(x)\; - \;
\bar\xi(V(x))\delta(V(x))L_{g}L_{G}V(x) L_{g}V(x)^\top.
\end{eqnarray}
We restrict our attention to functions $\delta$ such that
\begin{equation}
\label{300}
\delta'(V(x)) L_{G}V(x) \geq - \frac{1}{4} \; \; \forall x\in
{\mathbb R}^n.
\end{equation}
Recalling (\ref{2gh}) and (\ref{3h0}) therefore gives the
inequality
\begin{equation}
\label{bze}
\begin{array}{rcl}
\dot U & \leq & \frac{3}{4}\left[L_{f}V(x) - \bar\xi(V(x))
|L_{g}V(x)|^2\right]  \; +\;   \delta(V(x)) L_{f} L_{G}V(x)
\\
& &+ \bar\xi(V(x))\delta(V(x)) |L_{g}L_{G}V(x)||L_{g}V(x)|.
\end{array}
\end{equation}
>From (\ref{2gh}), we deduce that
\begin{eqnarray}
\label{5h0}
\dot U &\leq & \frac{1}{2}\left[L_{f}V(x) - \bar\xi(V(x))
|L_{g}V(x)|^2\right]  +    \delta(V(x)) L_{f} L_{G}V(x)\nonumber\\
&&+ \bar\xi(V(x))\delta^2(V(x)) |L_{g}L_{G}V(x)|^2 \ . \
\end{eqnarray}

\noindent {\em Third step.} The remaining part of the proof relies
extensively on the following:
\begin{lemma}
\label{tcl} Assume that the system (\ref{1gh}) satisfies
Assumptions  \ref{H1}-\ref{H2}. Then, there exist continuous
positive definite functions $\Gamma$ and $N$  satisfying the
following:  If $|L_{g}V(x)| \leq \Gamma(|x|)$, then either
$L_{f}V(x) \leq - N(|x|)$ or $L_{f}L_{G}V(x) \leq - N(|x|)$.
\end{lemma}\medskip
\begin{proof}
We first show that the continuous function \begin{equation}
\label{fb1}
\begin{array}{rcl}
 S(x) &=& \min\{0, L_{f}L_{G}V(x)\}+ \min\{0,
L_{f}V(x)\}- |L_{g}V(x)|
\end{array}\end{equation}
is negative definite. Observe first that $S(0) = 0$ and $S(x) \leq
0$ for all $x$. Assume  $S(x) = 0$. Each term of $S(x)$ is
nonpositive, so $\min\{0, L_{f}L_{G}V(x)\} = \min\{0, L_{f}V(x)\}
= |L_{g}V(x)| = 0$. By Assumption \ref{H2}, $x = 0$, which gives
the negative definiteness. Therefore $- S(x)$ is positive
definite, so we  can determine a continuous positive definite
real-valued function $\rho$ such that $\rho(|x|)
 \leq  - S(x)$ (e.g., $\rho(s)=\min\{-S(r): |r|=s\}$).
We prove that $|L_{g}V(x)| \leq \frac{1}{2}\rho(|x|)$ implies that
either $L_{f}L_{G}V(x) \leq - \frac{1}{4}\rho(|x|)$ or $L_{f}V(x)
\leq - \frac{1}{4}\rho(|x|)$. To this end, consider $x$ such that
$|L_{g}V(x)| \leq \frac{1}{2}\rho(|x|)$. Then  $\rho(|x|)  \leq -
\min\{0, L_{f}L_{G}V(x)\} - \min\{0, L_{f}V(x)\}  +
\frac{1}{2}\rho(|x|)$, by our choices of $\rho$ and $S$. We deduce
that $\min\{0, L_{f}L_{G}V(x)\} + \min\{0, L_{f}V(x)\} \; \leq \; -
\frac{1}{2}\rho(|x|)$. It follows that either $\min\{0,
L_{f}L_{G}V(x)\} \leq - \frac{1}{4}\rho(|x|)$ or $\min\{0, L_{f}
V(x)\} \leq - \frac{1}{4}\rho(|x|)$. Therefore, $|L_{g}V(x)| \leq
\frac{1}{2}\rho(|x|)$ implies $L_{f}L_{G}V(x) \leq -
\frac{1}{4}\rho(|x|)$ or $L_{f} V(x) \leq - \frac{1}{4}\rho(|x|)$,
so we can take $\Gamma(s)=\frac{1}{2}\rho(s)$ and
$N(s)=\frac{1}{4}\rho(s)$.
\end{proof}

\noindent {\em Fourth step.} We prove that the right hand side of
(\ref{5h0}) is negative definite when the smooth positive definite
function $\delta(\cdot)$ is suitably chosen.  By the preceding
lemma, there are three cases:

\noindent \underline{First Case.} $|L_{g}V(x)| \leq \Gamma(|x|)$
and $L_{f}V(x) \leq - N(|x|)$. Then the inequality (\ref{5h0})
implies that
\begin{equation}
\label{4h0}
\begin{array}{rcl}
\dot U  & \leq & - \frac{1}{2}N(|x|) + \delta(V(x)) L_{f}
L_{G}V(x) \; + \;   \bar\xi(V(x))\delta^2(V(x)) |L_{g}L_{G}V(x)|^2
\ .
\end{array}
\end{equation}
Choosing $\delta(\cdot)$ such that

\begin{equation}
\label{6h0} \delta(V(x)) L_{f} L_{G}V(x)  \leq
\frac{1}{8}N(|x|),\; \; \; \; \bar\xi(V(x))\delta^2(V(x))
|L_{g}L_{G}V(x)|^2 \le \frac{1}{8}N(|x|)
\end{equation}
for all $x\in {\mathbb R}^n$. Therefore, (\ref{4h0}) gives $\dot U
\leq - \frac{1}{4}N(|x|)< 0$ for all $x\neq 0$.

\noindent \underline{Second Case.} $|L_{g}V(x)| \leq \Gamma(|x|)$ and
$L_{f}L_{G}V(x) \leq - N(|x|)$. Then the inequalities (\ref{2gh})
and (\ref{5h0}) imply $\dot U   \leq  - \delta(V(x)) N(|x|)  +
\bar\xi(V(x))\delta^2(V(x)) |L_{g}L_{G}V(x)|^2$. Choosing
$\delta(\cdot)$ such that
\begin{equation}
\label{ns1} \delta(V(x)) \bar\xi(V(x))|L_{g}L_{G}V(x)|^2 \; \leq
\; \frac{1}{2}N(|x|)
\end{equation}
we obtain $\dot U  \leq - \frac{1}{2}\delta(V(x)) N(|x|)<0$ for
all $x\ne 0$.

\noindent \underline{Third Case.} $|L_{g}V(x)| \geq \Gamma(|x|)$.
Then the inequality (\ref{5h0}) implies that
\[
\begin{array}{rcl} \dot U  & \leq & - \frac{1}{2}\bar\xi(V(x))
\Gamma^2(|x|) \; + \; \delta(V(x)) L_{f} L_{G}V(x) \; + \;
\bar\xi(V(x))\delta^2(V(x)) |L_{g}L_{G}V(x)|^2 \ .
\end{array}
\]
Arguing as above provides  $\delta(\cdot)$ such that
\begin{equation}
\label{zs1}
\delta(V(x)) L_{f}
L_{G}V(x)  \leq  \frac{1}{8}\Gamma^2(|x|)\bar\xi(V(x)),\; \; \; \;
\delta^2(V(x)) |L_{g}L_{G}V(x)|^2   \leq \frac{1}{8}\Gamma^2(|x|)
\ ,
\end{equation}
so we obtain $\dot U \leq - \frac{1}{4}\bar\xi(V(x))
\Gamma^2(|x|) < 0$ for all $x \neq 0$.

\noindent {\em Fifth step.} To conclude the proof for the control
affine case, one has to prove that one can determine a $C^1$ and
positive definite function $\delta(\cdot)$ simultaneously satisfying
the requirements (\ref{7gh}), (\ref{300}), (\ref{6h0}), (\ref{ns1}),
(\ref{zs1}). This can be done as follows.  We can first find a $C^1$
positive definite function $\delta$ satisfying the requirements
(\ref{7gh}), (\ref{6h0}), (\ref{ns1}), (\ref{zs1}) that is
increasing on $[0,1]$, non-increasing on $[1,\infty)$ and bounded by
$1$. We denote this initial choice of $\delta$ by $\delta_a$. Next,
we minorize $1/(1+4|L_GV(x)|)$ by a positive function of the form
$x\mapsto {\cal P}(V(x))$ (using, e.g., ${\cal P}(s) =
\inf\{1/(1+4|L_GV(x)|): x\in {\mathbb R}^n, V(x)=s\}$). One can
easily determine an {}everywhere positive, non-increasing $C^1$
function $\omega(\cdot)$ such that $\omega(s) \leq
\frac{1}{2}\min\left\{{\cal P}(s), {\cal P}(2s), 1\right\}$ for all
$s \geq 0$. Now consider the function
\begin{equation}
\label{rtu1}
\delta(s) = \displaystyle\int_{\frac{1}{2}s}^{s} \frac{\delta_a(l)\delta_a(2l)\omega(l)}{1 + 4 l^2} dl \ .
\end{equation}
It is positive definite, of class $C^2$, and {}(since $\delta_a$ is
bounded by $1$) satisfies, for all $s \geq 0$,
\begin{equation}
\label{1tu1}
\begin{array}{rcl}
|\delta'(s)| & = &
\left|\frac{\delta_a(s)\delta_a(2s)\omega(s)}{1 + 4 s^2} - \frac{1}{2}\frac{\delta_a(\frac{1}{2}s)\delta_a(s)\omega(\frac{1}{2}s)}{1 + s^2}
\right|
\leq \omega(s) + \frac{1}{2} \omega(\frac{1}{2}s) \leq {\cal P}(s) \ .
\end{array}
\end{equation}
>From this inequality, one can deduce that $\delta$ defined in
(\ref{rtu1}) satisfies (\ref{300}). On the other hand, since
$\omega$ is smaller than $1$, the inequality {}$\delta(s) \leq
\int_{s/2}^{s} \delta_a(l)\delta_a(2l)/(1 + 4 l^2) dl$ is satisfied
for all $s \geq 0$.
Now, we distinguish between two cases. {\em First case:} {}If $s \in
[0,1]$, then, since $\delta_a$ is a nonnegative function smaller
than $1$ and increasing on $[0,1]$, we get {}$\delta(s) \leq
\int_{s/2}^{s} \delta_a(l) dl \leq \delta_a(s)$. {}{\em Second
case:} $s \geq 1$, then, since $\delta_a$ is a nonnegative function
smaller than $1$ and nonincreasing on $[1, + \infty)$, {}we get
\begin{equation}
\label{4tu1}
\delta(s) \leq \displaystyle\int_{\frac{1}{2}s}^{s} \frac{\delta_a(2l)}{1 + s^2} dl
\leq  \frac{s}{2(1+s^2)}\delta_a(s) \leq \delta_a(s)\ .
\end{equation}
Hence, the function $\delta$ defined in (\ref{rtu1}) satisfies
the requirements (\ref{7gh}), (\ref{300}), (\ref{6h0}), (\ref{ns1}), (\ref{zs1}).
\begin{remark}
\label{choose} The proof of Lemma  \ref{tcl} provides explicit
formulae for the functions $\Gamma$ and $N$ required for our
constructions. On the other hand, the function $\delta$ in
(\ref{3gh}) can be obtained by simply verifying the requirements in
the fifth step of our proof.
\end{remark}

\subsection{Fully Nonlinear Case}\label{fnc}
We now extend the construction to our original fully nonlinear
system (\ref{fn}). We can write
\begin{equation}
\label{expand}
F(x,u) = f(x) + g(x)u + h(x,u)u,\; \; {\rm where}\; \;
h(x,u)=\int_0^1\left[\frac{\partial F}{\partial u}(x,\lambda
u)-\frac{\partial F}{\partial u}(x,0)\right]{\rm d}\lambda.
\end{equation}
Along the trajectories of (\ref{fn}), it follows that $\dot{V} =
L_{f} V(x) + L_{g} V(x) u + \nabla V(x) h(x,u) u$. Since  $F$ is
$C^2$ in $u$, we can find a continuous function $H:[0,\infty)\times
[0,\infty)\to (0,\infty)$ that is nondecreasing in both variables
such that $|h(x,u) u| \le H(V(x),|u|)|u|^2$ for all $x$ and $u$. One
can find $\alpha_4 \in {\cal K}_{\infty}$ such that $|\nabla V(x)|
\leq \alpha_4(|x|)$ for all $x$. Taking $u$ to be a {}feedback of
the form (\ref{bj1}) gives
\[
%\label{8m}
%\begin{array}{rcl}
\dot{V} \; \leq\; L_{f} V(x) - \xi(V(x)) |L_{g}
V(x)|^2 + H_*(V(x),|\xi(V(x)) L_{g}V(x)^\top|)\xi^2(V(x))|L_{g} V(x)|^2
%\end{array}
\]
%\end{equation}
with $H_*(r,s) = \alpha_4(\alpha_1^{-1}(r)) H(r,s)$. We now restrict
our attention to the set ${\cal F}[\bar\xi]$ of all feedbacks
(\ref{bj1}) such that $\xi(s)  \leq \overline{\xi}(s)$ for all $s\ge
0$, where we {}assume  the positive function $\overline{\xi}$ is
such that
\begin{equation}
\label{10m}
H_*(V(x),\overline{\xi}(V(x))|L_{g}V(x)|)\overline{\xi}(V(x)) \;
\leq \; \frac{1}{2}\; \; \; \forall x\in {\mathbb R}^n.
\end{equation}
Condition (\ref{10m}) can  be satisfied by minorizing $\bar \xi$
as necessary without relabelling. (The proof that $\bar \xi$ can
be chosen to satisfy (\ref{10m}) is similar to the construction of
the function $\delta$ in the first part of the proof.) Fixing a
feedback from this family ${\cal F}[\bar\xi]$, we get
\begin{equation}
\label{11m} \dot{V} \leq L_{f} V(x) - \frac{1}{2} \xi(V(x)) |L_{g} V(x)|^2
\end{equation}
along the closed loop trajectories of (\ref{fn}). Applying the
construction from the first part of the proof  to the control
affine system (\ref{1gh}) with $\bar \xi=\xi$ provides
 a function $\delta$ and a CLF $U$ of the form (\ref{3gha}) such that
$W(x) := - \left\{L_{f} U(x) -  L_gU(x)\xi(V(x)) L_{g}V(x)^\top\right\}$
is positive definite.
Therefore, $\dot U$ along the trajectories of (\ref{fn}) in
closed-loop with the feedback (\ref{bj1}), with  $\xi$ satisfying
$\xi(s)\le \bar \xi(s)$ for all $s\ge 0$, reads
$\dot{U} = - W(x) - \nabla U(x)
h(x,- \xi(V(x)) L_{g}V(x)^\top) \xi(V(x)) L_{g}V(x)^\top$.
Therefore, since $H$ is non-decreasing in its second argument, it
follows from our choices of $\xi$ and $H$ that $\dot{U} \leq -
W(x) + \left\{|\nabla U(x)|
H(V(x),\overline{\xi}(V(x))|L_{g}V(x)|)\overline{\xi}(V(x))\right\}
\xi(V(x))\left|L_{g}V(x)\right|^2$ for all $x$. One can construct
a positive nondecreasing function $\Delta$ such that, along the
closed loop trajectories,
\begin{equation}
\label{15m}
\dot{U} \; \; \leq\; \;  - W(x) + \Delta(V(x))\xi(V(x))\left|L_{g}V(x)\right|^2.
\end{equation}
Now consider the function
(\ref{3gh}) with the above choice of $\delta$ and
$\Omega(s)=4\Delta(s)$, which is positive definite and radially
unbounded. Then, according to our Assumption \ref{H1},
(\ref{11m}), and (\ref{15m}),  we get
\[
\begin{array}{rcl}
\dot V^\sharp &\leq & - W(x) - \Delta(V(x))\xi(V(x)) |L_{g} V(x)|^2\; \;
 \forall x\in {\mathbb R}^n
\end{array}
\]The right-hand-side of this inequality is negative definite, so
we can satisfy the requirements of the theorem using
$\Omega(s)=4\Delta(s)$ and the CLF $V^\sharp$. This concludes our
proof.
\section{Example}
\label{secex}

We illustrate Theorem \ref{thgh} by applying it to the two-link
manipulator (see \cite{dang}). This system is a fully actuated
system described by the {}Euler-Lagrange equations
\begin{equation}
\label{ec1}
\begin{array}{l}
\left(m r^2 + M\frac{L^2}{3}\right)\ddot{\theta} + 2 M r
\dot{r}\dot{\theta} \;  = \;  \tau \ ,\; \; \; \;  m \ddot{r} - m r
\dot{\theta}^2 \;= \; F \ ,
\end{array}
\end{equation}
where $M$ is the mass of the arm; $L$ is its length; $m$ is the mass
of the gripper; $r$ and $\theta$ denote  the angle of the link and
the position of the gripper, respectively;  and $\tau$ and $F$ are
forces acting on the system. It is well-known that (\ref{ec1}) can
be stabilized by bounded control laws. On the other hand, this
system is globally feedback linearizable so  a quadratic CLF can be
determined. The novelty  is that we determine a CLF whose derivative
along the trajectory is made negative definite by an appropriate
choice of {\em bounded} feedback. Without loss of generality, we
take $m = M = 1$ and $L = \sqrt{3}$. With $x_1 := \theta, x_2 :=
\dot{\theta}, x_3:= r, x_4 := \dot{r}$, the system (\ref{ec1})
becomes
\begin{equation}
\label{msa1}
\begin{array}{l}
\dot{x}_1  =  x_2,\; \; \;  \; \dot{x}_2  =  - \frac{2 x_3 x_2
x_4}{x_3^2 + 1} + \frac{\tau}{x_3^2 + 1},\; \;  \; \; \dot{x}_3  =
x_4, \; \; \; \; \dot{x}_4 =  x_3 x_2^2 + F.
\end{array}
\end{equation}
We construct a globally stabilizing feedback, bounded in norm by
$1$, and an associated CLF for (\ref{msa1}).  We set $\langle
p\rangle =1/(2\sqrt{1+p^2}\, )$ for all $p\in {\mathbb R}$
throughout the sequel.

Consider the function
\begin{equation}
\label{msa2}
V(x) = \frac{1}{2}\left[(x_3^2 + 1)x_2^2 + x_4^2 + \sqrt{1 + x_1^2} + \sqrt{1 + x_3^2} - 2\right] \ .
\end{equation}
This function is composed of the kinetic energy of the system with
additional terms. It is positive definite and radially unbounded and
its derivative along trajectories of (\ref{msa1}) satisfies
$\dot{V}(x) = x_2 \tau + x_4 F + x_1 \langle x_1 \rangle x_2 +x_3
\langle x_3 \rangle x_4$.
 Therefore the change of feedback
\begin{equation}
\label{msa4} \tau = - x_1\langle x_1\rangle + \tau_b \; , \; F =
-x_3\langle x_3\rangle  + F_b
\end{equation}
yields $\dot{V}(x) = x_2 \tau_b + x_4 F_b$.
On the other hand, after the change of feedback (\ref{msa4}), the
equations of the system take the control affine form $\dot
x=f(x)+g(x)u$ with
\[
f(x)= \left[
\begin{array}{l}
x_2\\
\frac{- 2 x_3 x_2 x_4 - x_1\langle x_1\rangle}{x_3^2 + 1}\\
x_4\\
x^2_2x_3-x_3\langle x_3\rangle
\end{array}\right],\; \; \; \;
g(x)=\left[ \begin{array}{ll} 0 &0\\
\frac{1}{{}x^2_3+1}& 0\\
0& 0\\
0& 1
\end{array}\right],\; \; \; \; u=\left[\begin{array}{l}\tau_b\\
F_b\end{array}\right]
\]
Next consider the vector field $G(x) = (0, x_1, 0, x_3)^\top$.
Simple calculations yield
\begin{equation}
\label{msa8} L_G V(x) = \frac{\partial V}{\partial x_2}(x) x_1 +
\frac{\partial V}{\partial x_4}(x) x_3 = (x_3^2 + 1) x_2 x_1 + x_4
x_3.
\end{equation}
Since $\nabla (L_G V(x))=(x_2(x^2_3+1), x_1(x^2_3+1),
x_4+2x_1x_2x_3, x_3)$, we get
\begin{equation}
\label{msa9}
\begin{array}{rcl}
L_f L_G V(x) & = & x_2^2 (2 x_3^2 + 1) + x_4^2 - x_1^2 \langle x_1\rangle
- x_3^2 \langle x_3\rangle .
\end{array}
\end{equation}
We now check that Assumptions \ref{H1} and \ref{H2} are satisfied.
Since $L_f V(x) = 0$ and  $L_g V(x) = [x_2 \;\; x_4]$ everywhere,
Assumption \ref{H1} is satisfied. If $L_{g}V(x) = 0$, then $x_2 =
x_4 = 0$, in which case we get $L_f L_G V(x) = - x^2_1\langle
x_1\rangle - x^2_3\langle x_3\rangle$. It follows that if $x \neq 0$
and $L_{g}V(x) = 0$, then $L_f L_G V(x) < 0$. Therefore Assumption
\ref{H2} is satisfied. {}Hence, Theorem \ref{thgh} applies. Consider
the function
\begin{equation}
\label{dra1} V^\sharp(x) \; = \; 40[2 + 2V(x)]^6 + L_{G}V(x) -
{}40(2^6) \ .
\end{equation}
Simple multiplications show  $80[2 + 2V(x)]^6  \geq V^\sharp(x) \geq
3\left(x_1^2 + x_2^2 + x_3^2 + x_4^2\right)$ for all $x$, so
$V^\sharp$  is positive definite and radially unbounded. Moreover,
 we see that along the
trajectories of (\ref{msa1}) after the change of feedback
(\ref{msa4}),
\begin{equation}
\label{3jlk}
\begin{array}{rcl}
\dot{V}^\sharp(x) & = & 480[2 + 2V(x)]^5 (x_2 \tau_b + x_4 F_b) + x_2^2 (2 x_3^2 + 1) + x_4^2
- x^2_1\langle x_1\rangle - x^2_3\langle x_3\rangle
\\
& & + x_1 \tau_b + x_3 F_b \ ,
\end{array}
\end{equation}
since  $\dot{V}(x) = x_2 \tau_b + x_4 F_b$.  {}Therefore, from the
triangle inequality, we deduce that
\begin{equation}
\label{tri1}
\begin{array}{rcl}
\dot{V}^\sharp(x) & \leq & \sqrt{1 + x_1^2} \tau_b^2 +
480[2 + 2V(x)]^5 x_2 \tau_b +  x_2^2 (2 x_3^2 + 1)
\\
& & + \sqrt{1 + x_3^2} F_b^2 + 480[2 + 2V(x)]^5 x_4 F_b + x_4^2
- \frac{1}{2}x^2_1\langle x_1\rangle - \frac{1}{2}x^2_3\langle x_3\rangle \ .
\end{array}
\end{equation}
We demonstrate now that $V^\sharp$ is a CLF for (\ref{msa1}) by
showing that the right hand side of {}(\ref{tri1}) is negative
definite for the feedbacks
\begin{equation}
\label{ldr1} \tau_b = - x_2\langle x_2\rangle \; , \; F_b = - x_4\langle x_4\rangle.
\end{equation}
To this end, notice that we have
\begin{equation}
\label{3khk}
\begin{array}{rcl}
\dot{V}^\sharp(x) & \leq & T_1(x) x_2^2\langle x_2\rangle + T_2(x) x_4^2\langle x_4\rangle
- \frac{1}{2}\left[x_1^2\langle x_1\rangle + x_2^2\langle x_2\rangle + x_3^2\langle x_3\rangle + x_4^2\langle x_4\rangle\right]
\end{array}
\end{equation}
where we define the $T_i$'s by $T_1(x) = \sqrt{1 + x_1^2} - 480(2 +
2 V(x))^5 + 2 \sqrt{1 + x_2^2}(2 x_3^2 + 1) + \frac{1}{2}$ and
$T_2(x) = \sqrt{1 + x_3^2} - 480(2 + 2 V(x))^5 + 2 \sqrt{1 + x_4^2}
+ \frac{1}{2}$. From the expression of $V(x)$, we deduce that $T_1$
and $T_2$ are nonpositive and therefore
\begin{equation}
\label{3jhk}
\dot{V}^\sharp(x) \leq - \frac{1}{2}\left[x_1^2\langle x_1\rangle + x_2^2\langle x_2\rangle + x_3^2\langle x_3\rangle + x_4^2\langle x_4\rangle \right]\ .
\end{equation}
The right hand side of this inequality is negative definite and the feedbacks resulting from (\ref{msa4})
and (\ref{ldr1}) are bounded in norm by $1$. This concludes the proof.

\section{Robustness to Actuator Errors}
\label{seciss}
Theorem \ref{thgh} provided a stabilizing
feedback  $u=K_1(x)$  such that $\dot x=f(x)+g(x)K_1(x)$ is
globally asymptotically stable (GAS) to $x=0$. Moreover, for each
$\varepsilon>0$, we can choose $K_1$ to satisfy $|K_1(x)|\le
\varepsilon$ for all $x\in {\mathbb R}^n$.

One natural and widely used generalization of the GAS condition is
the so-called input-to-state stable (ISS) property \cite{S89}. For
a general nonlinear system $\dot x=F(x,d)$ evolving on ${\mathbb
R}^n\times {\mathbb R}^m$ (where $d$ represents the disturbance),
the ISS property is the requirement that there exist $\beta\in
{\cal KL}$ and $\gamma\in {\cal K}_\infty$ such that the following
holds for all measurable essentially bounded functions ${\mathbf
d}:[0,\infty)\to{\mathbb R}^m$ and corresponding trajectories
$x(t)$ for $\dot x(t)=F(x(t),{\mathbf d}(t))$:
\begin{equation}
\tag{ISS} |x(t)|\le \beta(|x(0)|,t)+\gamma(|{\mathbf d}|_\infty)\;
\; \forall t\ge 0.
\end{equation}
Here $|\cdot|_\infty$ is the essential supremum norm.
The ISS property reduces to GAS to $0$ for systems with no
controls, in which case the overshoot term $\gamma(|{\mathbf
d}|_\infty)$ in the ISS decay condition  is $0$;  see also
\cite{MRS04, MS04} for the relationship between the ISS property
and asymptotic controllability.  It is therefore natural to look
for a feedback $K(x)$ for (\ref{1gh})
(which could in principle differ from  $K_1$) for which
\begin{equation}
\label{open} \dot x=F(x,d):=f(x)+g(x)[K(x)+d]
\end{equation}
is ISS, and for which $|K(x)|\le \varepsilon$ for all $x\in
{\mathbb R}^n$, where $\varepsilon$ is any prescribed positive
constant. In other words, we would want an arbitrarily small
feedback $K$ that renders (\ref{1gh}) GAS to $x=0$ and that has
the additional property that (\ref{open}) is also ISS with respect
to actuator errors $d$.

However, it is clear that this objective cannot be met, since
there is no {\em bounded} feedback $K(x)$ such that the
one-dimensional system $\dot x=K(x)+d$ is ISS. On the other hand,
if we add

\begin{assumption}
\label{H3}
 A positive
nondecreasing smooth function ${\cal D}$ such that (i)
$\int_{0}^{+\infty} \frac{1}{{\cal D}(s)}\, {\rm d} s = +\infty$
and
 (ii) $|L_gV(x)|\le {\cal D}(V(x))$ for
all $x\in {\mathbb R}^n$ is known.
\end{assumption}

\noindent where $V$ satisfies our continuing Assumptions
\ref{H1}-\ref{H2}, then any feedback $K:=-\xi(V(x))L_gV(x)^\top$,
obtained from Theorem \ref{thgh} for the control affine system
$\dot x=f(x)+g(x)u$ and chosen such that $|\xi(V(x))L_gV(x)|\le
\varepsilon$ for all $x\in {\mathbb R}^n$, also renders
(\ref{open}) {\em integral-input-to-state stable (iISS)}.
%(This
%fact is easy to check when $V$ is a CLF for the system but extra
%care is needed because $V$ might not be a CLF.)
 For a general
nonlinear system $\dot x=F(x,d)$ evolving on ${\mathbb R}^n\times
{\mathbb R}^m$, the iISS condition is the following: There exist
$\beta\in {\cal KL}$ and $\alpha,\gamma\in {\cal K}_\infty$ such
that for all measurable locally  essentially bounded functions
${\mathbf d}:[0,\infty)\to{\mathbb R}^m$ and corresponding
trajectories $x(t)$ for $\dot x(t)=F(x(t),{\mathbf d}(t))$,
\begin{equation}
\tag{iISS} \alpha(|x(t)|)\le
\beta(|x(0)|,t)+\int_0^t\gamma(|{\mathbf d}(s)|){\rm d}s\; \; \;
\;  \forall t\ge 0.
\end{equation}
The iISS condition reflects the qualitative property of having
small overshoots when the disturbances have finite energy. It
provides a nonlinear analog of ``finite $H^2$ norm'' for linear
systems, and thus has obvious physical relevance and significance
\cite{AAS02,  ASW00}.
Assumptions \ref{H1}-\ref{H3} hold for our example
in the previous section, since in that case, $|L_gV(x)|\le
2(V(x)+2)$ for all $x\in {\mathbb R}^n$, so we can take ${\cal
D}(s)=2(s+2)$.  In fact, our assumptions hold for a broader class
of Hamiltonian systems as well;  see Remark \ref{newrk} below.

To verify that the Theorem \ref{thgh}  feedback also renders
(\ref{open}) iISS, we begin by fixing $\varepsilon>0$ and $V$
satisfying our Assumptions \ref{H1}-\ref{H3}, and applying our
theorem to $\dot x=f(x)+g(x)u$. This provides a CLF $U$ for
(\ref{1gh}) and  a corresponding positive function $\xi$ that
satisfies $|\xi(V(x))L_gV(x)|\le \varepsilon$ for all  $x\in
{\mathbb R}^n$.  The CLF $U$ has the form (\ref{3gha}). By
reducing $\delta$ and $\delta'$ from Section \ref{ca}, and
replacing ${\cal D}(p)$ with $p\mapsto {\cal D}(2p)+1$ in
Assumption \ref{H3} without relabelling, we can assume
\begin{equation}
\label{key} |L_gU(x)|\le {\cal D}(U(x)) \; \; \forall x
\in {\mathbb R}^n.\end{equation} Then
 \begin{equation}\label{newclf}
 \tilde
U(x)=\int_0^{U(x)}\frac{ {\rm d}p}{{\cal D}(p)},\; \; {\rm
where}\; \;  U(x)=V(x)+\delta(V(x))L_GV(x)
\end{equation}
 is again a
CLF for our dynamic (\ref{1gh}), since  our choice of ${\cal D}$
gives $\tilde U(x)\to +\infty$ as $|x|\to\infty$ because $U$ is
radially unbounded, and  because $\nabla \tilde U(x)\equiv \nabla
U(x)/{\cal D}(U(x))$ (which gives the CLF decay condition).
The smoothness of $\tilde U$ follows because $U$ and ${\cal D}$ are
both smooth. Finally, (\ref{key}) gives
\begin{equation}
\label{key3}
|L_g\tilde U(x)|=|L_gU(x)/{\cal D}(U(x))|\le 1\; \; \forall x
\in {\mathbb R}^n.
\end{equation}
We next choose the smooth feedback
$K_1(x) = - \xi(V(x))L_gV(x)^\top$,
where $\xi$ is a smooth positive function satisfying the above
requirements, so $K_1$ renders (\ref{1gh}) GAS to $x=0$, by
Theorem \ref{thgh}. To check that $K(x):=K_1(x)$ also renders
(\ref{open}) iISS, notice that our choice of $K_1$ and
(\ref{key3}) give
\begin{equation}
\label{decay}
\begin{array}{rcl} \nabla \tilde
U(x) F(x,d) &=& \nabla \tilde U(x) [f(x)+g(x)K_1(x)]+ L_g\tilde
U(x)d \\&\le& -\alpha_5(|x|)+|L_g\tilde U(x)|\, |d|\; \;  \le\; \;
-\alpha_5(|x|)+\, |d|
\end{array}
\end{equation}
for all $x$ and $d$ and some continuous positive definite function $\alpha_5$.
Inequality (\ref{decay}) says (see \cite{ASW00}) that the positive definite
radially unbounded smooth function $\tilde U$ is an iISS-CLF for (\ref{open}).
The fact that (\ref{open}) is iISS  now follows from the iISS
Lyapunov characterization \cite[Theorem 1]{ASW00}.  We conclude as
follows:
\begin{corollary}
Let the data (\ref{data}) satisfy Assumptions \ref{H1}-\ref{H3}
for some vector field $G:{\mathbb R}^n\to {\mathbb R}^n$ and
$V:{\mathbb R}^n\to {\mathbb R}$, and let $\varepsilon>0$ be
given. Then there exist smooth functions $\delta,\xi:[0,\infty)\to
[0,\infty)$ such that  (i) the system
 (\ref{open})   with  the feedback $K(x):=-\xi(V(x))L_gV(x)^\top$ is iISS and has a smooth iISS-CLF
of the form (\ref{newclf}) and  (ii) $|K(x)|\le \varepsilon$ for
all $x\in {\mathbb R}^n$.
\end{corollary}

\begin{remark}
\label{newrk} Assume the Hamiltonian system (\ref{hamil})
satisfies the conditions (a)-(b) we introduced in Section
\ref{discussion} as well as the following additional condition:
(c) There exist $\underline \lambda, \bar\lambda>0$ such that
${\rm spectrum}\{M^{-1}(q)\}\subseteq [\underline\lambda, \bar
\lambda]$ for all $q$.  (Assumption (c) means there are {\em
positive} constants $\underline c$ and $\bar c$ such that
$\underline c|p|^2\le p^\top M(q)p\le \bar c|p|^2$ for all $q$ and
$p$.) Then (\ref{hamil}) satisfies our Assumptions
\ref{H1}-\ref{H3} and so is covered by the preceding corollary.
 In fact, we saw on p. \pageref{hamil} that (a)-(b) imply that Assumptions
\ref{H1}-\ref{H2} hold with $V=H$, and then Assumption \ref{H3}
follows from (c) because $|L_gV(x)|^2=|p^\top M^{-1}(q)|^2\le \bar
\lambda^2 |p|^2\le (\bar \lambda^2/\underline \lambda)p^\top
M^{-1}(q)p\le 2 (\bar \lambda^2/\underline \lambda) V(x)$ for all
$x=(q,p)$. We can choose
${\cal D}(s):=\sqrt{2 (\bar \lambda^2/\underline\lambda)(s+1)}$.
\end{remark}

\section{Conclusion}
\label{seconc}

We showed how to construct control-Lyapunov functions for fully
nonlinear systems satisfying  appropriate generalizations of the
Jurdjevic-Quinn conditions. We also constructed
 feedbacks  of arbitrarily small norm that
render our systems integral-input-to-state stable to actuator
errors. Our constructions apply to important families of nonlinear
systems, and in particular to systems described by Euler-Lagrange
equations. Redesign and further robustness analysis for our
systems via our construction of control-Lyapunov functions will be
subjects of  future work.

%%%%%%%%%%%%%%%%%%%%%%%%%%%%%%%%%%%%%%%%%%%%%%%%%%%%%%%%%%%%%%%%%%%%%%%%%%%%%%%%
\section{Acknowledgements}

The authors thank J.-B. Pomet for useful discussions.

%%%%%%%%%%%%%%%%%%%%%%%%%%%%%%%%%%%%%%%%%%%%%%%%%%%%%%%%%%%%%%%%%%%%%%%%%%%%%%%%

\end{document}